\documentclass[a4paper]{amsproc}
\usepackage{amsfonts}
\usepackage{amssymb}
\usepackage{amscd}
\usepackage{verbatim}
\usepackage{graphicx}
\usepackage{subfigure}
\usepackage{fancyhdr}
\usepackage[T1]{fontenc}
\usepackage[latin1]{inputenc}
\usepackage[english]{babel}

\input xy
\xyoption{all}

\newcommand{\fait}[3]{\begin{#1}\label{#2}{#3}\end{#1}}
\newcommand{\cqfd}{\hfill \mbox{$\square$}}

\def\Hom{\mathop{\rm Hom}\nolimits}

\def\End{\mathop{\rm End}\nolimits}
\def\Aut{\mathop{\rm Aut}\nolimits}

\def\Ext{\mathop{\rm Ext}\nolimits}

\def\path{\mathop{\rightsquigarrow}\nolimits}
\def\mod{\mathop{\rm mod}\nolimits}

\def\add{\mathop{\rm add}\nolimits}

\def\proj{\mathop{\rm proj}\nolimits}

\newcommand{\R}{A[G]}

\newcommand{\cohX}{\mbox{coh}\mathbb{X}}

\newcommand{\DBA}{D^b(A)}
\newcommand{\DBB}{D^b(B)}
\newcommand{\DBH}{D^b(H)}

\newcommand{\DDBH}{D^b(\mc{H})}

\newcommand{\B}{\End_A T}
\newcommand{\Ga}{\Gamma}

\newcommand{\mc}{\mathcal}

\newcommand{\A}{\mathcal{A}}

\newcommand{\M}{{{M}^{\bullet}}}
\newcommand{\N}{{{N}^{\bullet}}}
\newcommand{\X}{{{X}^{\bullet}}}

\newcommand{\T}{{{T}^{\bullet}}}

\newcommand{\mor}[3]{\xymatrix@1@C=15pt{#1: #2\ar[r]& #3}}
\newcommand{\map}[2]{\xymatrix@1@C=15pt{#1\ar@{|->}[r]& #2}}

\begin{document}

\title[\textsl{Skew group algebras of piecewise hereditary algebras}]
{Skew group algebras of piecewise hereditary algebras are piecewise hereditary}%
\date{}
\author[J.~Dionne, M.~Lanzilotta and D.~Smith]{Julie Dionne, Marcelo Lanzilotta and David Smith }

\begin{abstract}
\noindent We show that the
main results of Happel-Rickard-Schofield (1988) and
Happel-Reiten-Smal{\o} (1996) on piecewise hereditary algebras are
coherent with the notion of group action on an algebra. Then, we
take advantage of this compatibility and
show that if $G$ is a finite group acting on a piecewise hereditary
algebra $A$ over an algebraically closed field whose characteristic
does not divide the order of $G$, then the resulting skew group
algebra $A[G]$ is also piecewise hereditary.
\end{abstract}

\bibliographystyle{plain}

\maketitle

%-------------------------------------------------------------------------------
%
%  SECTION  : INTRODUCTION
%
%-------------------------------------------------------------------------------
%\section{Introduction}
%    \label{Introduction}

Let $k$ be an algebraically closed field.  For a finite dimensional
$k$-algebra $A$, we denote by $\mod A$ the category of finite
dimensional left $A$-modules, and by $\DBA$ the (triangulated)
derived category of bounded complexes over $\mod A$. Let $\mc{H}$ be a connected hereditary abelian
$k$-category. Following
\cite{HRS96II} (compare \cite{H88,HRS88}), we say that $A$ is
\textbf{piecewise hereditary of type $\mc{H}$} if it is derived
equivalent to $\mc{H}$, that is $\DBA$ is triangle-equivalent to the
derived category $D^b(\mc{H})$ of bounded complexes over $\mc{H}$.
Over the years, piecewise hereditary algebras have been widely
investigated and
proved to be related with many other topics, such as the simply
connected algebras and the trivial extensions, the self-injective algebras of
polynomial growth and the strong global
dimension.

Hereditary categories $\mc{H}$ having tilting objects are of
special interest in representation theory of algebras.  The 
endomorphism algebras $\End_{\mc{H}}T$ of tilting objects $T$
in $\mc{H}$, called \textbf{quasitilted} algebras, were introduced
and studied in \cite{HRS96}.  It is well-known that $\mc{H}$ and
$\End_{\mc{H}}T$ are derived equivalent. When $k$ is algebraically closed, it was shown by Happel
\cite{Hap01} that $\mc{H}$ is either derived equivalent to a
finite dimensional hereditary $k$-algebra $H$ or derived
equivalent a category of coherent sheaves $\cohX$ on a weighted
projective line $\mathbb{X}$ (in the sense of \cite{GL87}).

%On the other hand, an \textbf{action} of a group $G$ on an algebra
%$A$ is a group homomorphism from $G$ to the group of automorphisms
%of $A$. For such an action, one defines the skew group algebra
%$A[G]$ (see Section \ref{Skew group algebras}). The study of the
%representation theory of skew group algebras was started in
%\cite{dlP83,RR85}, and pursued in \cite{ALR07,FR02,Smith07}. 

The aim of this paper is to study the skew group algebra $A[G]$ (see Section \ref{Skew group algebras}), in case $A$ is a piecewise
hereditary algebra. Our main result (Theorem \ref{thm DLS3}) shows
that under standard assumptions, the skew group algebra $A[G]$ is also
piecewise hereditary.

In order to give a clear statement of our main results, we need
additional terminology.
Let $G$ be a group and $\A$ be an additive category.  An
\textbf{action} of $G$ on $\A$ is a group homomorphism
$\mor{\theta}{G}{\Aut \A} \ (\map{\sigma}{\theta_\sigma})$ from $G$
to the group of automorphisms of $\A$.
An object $\mathcal M$ in $\A$ is \textbf{$G$-stable} with
respect to $\theta$, or briefly \textbf{$G$-stable} in case there
is no ambiguity, if $\theta_\sigma\mathcal M\cong\mathcal M$ for
all $\sigma\in G$. For such an object, the algebra
$B=\End_{\A}\mathcal M$ inherits an action of $G$ from $\theta$,
denoted $\theta_B$; see (\ref{rem endo}). Also, given another
additive category $\mc{B}$ and an action
$\mor{\vartheta}{G}{\Aut \mc{B}} \
(\map{\sigma}{\vartheta_\sigma})$
of $G$ on $\mc{B}$, a functor $\mor{F}{\A}{\mc{B}}$ is \textbf{$G$-compatible}, with respect to the pair $(\theta,
\vartheta)$, if $F\theta_\sigma=\vartheta_\sigma F$ for every
$\sigma\in G$.

Examples of particular interest occur when a group $G$ acts on an
artin algebra $A$ as above. Then the action $\theta$ of $G$ on $A$
induces an action $\theta_{\mod A}$ of $G$ on $\mod A$, and further an action $\theta_{D^b(A)}$ on
$D^b(A)$; see Sections \ref{Skew group algebras} and \ref{Group
actions on homotopy and derived categories}.
%

%Let $A$ and $B$ be algebras. In the vein of \cite{H88}, we say that
%$B$ is \textbf{(co)tiltable} to $A$ if there is a sequence $B=A_0,
%A_1, \dots, A_n=A$ of algebras and a sequence $T_0, T_1, \dots,
%T_{n-1}$ where $T_i$ is a tilting module or a cotilting module for
%$A_i$ with endomorphism ring isomorphic to $A_{i+1}$. If each $T_i$
%is a tilting module for $A_i$, we say that $B$ is \textbf{tiltable}
%to $A$. If $B$ is hereditary and we can further ensure that each
%$T_i$ is \textbf{splitting}, that is it induces a splitting torsion
%theory on $\mod A_{i+1}$, then we say that $A$ is \textbf{iterated
%tilted} of type $B$. In addition, if $\vartheta$ is an action of $G$
%on $B$ and $T_i$ is $G$-stable (with respect to the action of $G$
%on $A_i$ inherited from $\vartheta$) for each $i$, then we use the
%terminology \textbf{$G$-(co)tiltable}, \textbf{$G$-tiltable} or
%\textbf{$G$-iterated tilted algebra of type $B$} depending on the
%situation.

This opens the way to our main results. Our first theorem stands as a generalization of the main result in \cite{HRS88} and \cite{HRS96II}.

%*********** Theorem : thm DLS2
\fait{my_thm}{thm DLS2}{Let $A$ be a $k$-algebra, and 
$\mc{H}=\emph{mod}H$, with $H$ a hereditary algebra, or
$\mc{H}=\emph{coh}\mathbb{X}$, the category of coherent sheaves
on a weighted projective line $\mathbb{X}$. Let $G$ be a group, and 
$\mor{\theta}{G}{\Aut A}$ and $\mor{\vartheta}{G}{\Aut \mc{H}}$ be
fixed actions of $G$ on $A$ and $\mc{H}$.%
\begin{enumerate}
\item[\emph{(a)}] The following conditions are equivalent :
\begin{enumerate}
\item[\emph{(i)}] There exists a $G$-compatible triangle-equivalence
$\mor{E}{D^b(\mc{H})}{\DBA}$ (with respect to the pair of induced
actions $(\vartheta_{D^b(\mc{H})}, \theta_{D^b(A)})$);
\item[\emph{(ii)}] There exist a $G$-stable tilting object $T$ in $\mc{H}$ and sequences
$\End_{\mc{H}}T=A_0, A_1, \dots, A_n=A$ of $k$-algebras and $T_0,
T_1, \dots, T_{n-1}$ of modules such that, for each $i$,
$A_{i+1}=\End_{A_i}T_i$ and $T_i$ is a $G$-stable tilting or
cotilting $A_i$-module (with respect to the induced action
$\vartheta_{\mod A_i}$), and the induced action $\vartheta_{A_n}$
coincides with $\theta$;
\item[\emph{(iii)}] There exist a $G$-stable tilting object $T$ in $\mc{H}$ and sequences
$\End_{\mc{H}}T=A_0, A_1, \dots, A_n=A$ of $k$-algebras and $T_0,
T_1, \dots, T_{n-1}$ of modules such that, for each $i$,
$A_{i+1}=\End_{A_i}T_i$ and $T_i$ is a $G$-stable splitting
tilting or cotilting $A_i$-module (with respect to the induced
action $\vartheta_{\mod A_i}$), and the induced action
$\vartheta_{A_n}$ coincides with $\theta$;
\end{enumerate}
\item[\emph{(b)}] If the above conditions are satisfied and $G$ is a finite group whose order is
not a multiple of the characteristic of $k$, then the algebra
$(\End_{\mc{H}}T)[G]$, where $T$ is as in (ii) or in (iii), is
quasitilted and derived equivalent to $A[G]$. In particular,
$A[G]$ is piecewise
hereditary. %
\end{enumerate}}

The equivalence of the conditions (a)(i)-(a)(iii) of Theorem \ref{thm DLS2} were previously shown in \cite{HRS96II,HRS88} in the case where, essentially,
the actions $\theta$ and $\vartheta$ are the trivial actions, that
is trivial homomorphisms of groups. Actually, our proofs are
adaptations of the original ones.

In addition, it will become clear in Section \ref{Proof of Theorem
3} that any
triangle-equivalence $\xymatrix@1@C=15pt{D^b(\mc{H})\ar[r]& \DBA}$ can be converted
into a $G$-compatible triangle-equivalence.  As an application of
this observation, together with Theorem
\ref{thm DLS2} and Happel's Theorem \cite{Hap01}, we will obtain our
main theorem.

%*********** Theorem : thm DLS3
\fait{my_thm}{thm DLS3}{Let $A$ be a piecewise hereditary
$k$-algebra of type $\mathcal{H}$, for some Ext-finite hereditary
abelian $k$-category with tilting objects $\mc{H}$.  Moreover, let
$G$ be a finite group whose order is not a multiple of the
characteristic of $k$.  Then,
\begin{enumerate}
\item[\emph{(a)}] If $\mathcal{H}=\mod H$, for some hereditary algebra $H$, then for
any action of $G$ on $A$, there exist a hereditary algebra $H'$, derived equivalent to $H$,
and an action of $G$ on $H'$ such that $A[G]$ is piecewise hereditary of type $\mod
H'[G]$.
\item[\emph{(b)}] If $\mathcal{H}=\emph{coh}\mathbb{X}$, for some category of coherent sheaves on a weighted
projective line $\mathbb{X}$, then for any action of $G$ on $A$
there exist an action of $G$ on $\mathcal{H}$ and a $G$-compatible
triangle-equivalence $\mor{E}{D^b(\mc{H})}{\DBA}$.  In particular,
$A[G]$ is piecewise hereditary.%
\end{enumerate}}

%As we shall see, the algebras $H$ and $H'$ in the statement (a)
%above generally differ from each other.

In Section \ref{Preliminaries},
we fix the notations and terminologies. Most of the necessary background
on weighted projective lines is however postponed to Section
\ref{Proof of Theorem 3}, since it is not explicitly needed until
then. In Section \ref{Group actions and $G$-compatible derived
equivalences}, we study the $G$-compatible
triangle-equivalences of derived categories induced by the
$G$-stable tilting modules. Section \ref{Piecewise hereditary
algebras revisited} is devoted to the proof of Theorem \ref{thm DLS2}.  In Section \ref{Proof of Theorem
3}, we prove Theorem \ref{thm DLS3}. This involves showing
that any triangle-equivalence between $D^b(\mc{H})$ and $\DBA$ can
be converted into a $G$-compatible equivalence when $\mc{H}$ is a
module category over a hereditary algebra or a category of coherent
sheaves on a weighted projective line. Finally, in Section \ref{An
example}, we give an illustrative example.

%-------------------------------------------------------------------------------
%
%  SECTION 1 : PRELIMINARIES
%
%-------------------------------------------------------------------------------
\section{Preliminaries}
    \label{Preliminaries}

In this paper, all considered algebras are finite dimensional
algebras over an algebraically closed field $k$ (and, unless
otherwise specified, basic and connected). Moreover, all modules
are finitely generated left modules.
For an algebra $A$, we denote by $\proj A$ a full subcategory of $\mod
A$ consisting of one representative from each isomorphism class of
indecomposable projective modules. Given an $A$-module $T$, we let
$\add T$ be the full subcategory of $\mod A$ having as objects the
direct sums of indecomposable direct summands of $T$. Also, the
functor $D = \Hom_k(-, k)$ is the standard duality between $\mod
A$ and $\mod A^{op}$.

%-------------------------------------------------------------------------------
%  SUBSECTION 1.1 : TILTING THEORY
%-------------------------------------------------------------------------------
%\subsection{Tilting theory ? changer titre}
%    \label{Tilting theory}

Let $A$ be an algebra.  An $A$-module $T$ is a \textbf{tilting
module} if $T$ has projective dimension at most one,
$\Ext^1_A(T,T)=0$ and there exists a short exact sequence of
$A$-modules
$\xymatrix@1@C=15pt{0\ar[r]& A\ar[r] & T_0 \ar[r] & T_1\ar[r] &
0}$
in $\mod A$, with $T_0, T_1\in\add T$.

For basic results on tilting theorey, we refer to \cite{ASS06}, and for derived categories we refer to \cite{H88} or \cite{V77}. For an object $M$ in a triangulated category, we shall denote the image of $M$ under the "shift" self-equivalence $T$ by $M[1]$, and similarly $T^nM$ will be denoted by $M[n]$ for any $n$.

%In most of the situations, $\A$ will denote the category $\mod A$
%for some finite dimensional $k$-algebra $A$.  In this case, we
%simply denote $D^b(\mod A)$ by $D^b(A)$. The bounded derived
%category is best understood when $\A$ is a hereditary abelian
%category, that is such that the bifunctor $\Ext^2(-,?)$ vanishes.
%Then, any indecomposable object $\M=(M^i,
%d^i_{\M})_{i\in\mathbb{Z}}$ in $D^b(\A)$ is isomorphic to a stalk
%complex, and so we can assume that there exists $i_0\in\mathbb{Z}$
%such that $M^i\neq 0$ if and only if $i=i_0$, and where $M^{i_0}$ is
%an indecomposable object in $\A$. In this case, we write
%$\M=M^{i_0}[-i_0]$. Moreover, given two objects $M, N$ in $\A$ and
%$i\in \mathbb{Z}$, we have $\Hom_{D^b(\A)}(M[i], N[i])\cong
%\Hom_{\A}(M, N)$. This relation shows that derived categories of
%piecewise hereditary algebras behave quite well.

%-------------------------------------------------------------------------------
%  SUBSECTION 1.2 : SKEW GROUP ALGEBRAS
%-------------------------------------------------------------------------------
\subsection{Skew group algebras}
    \label{Skew group algebras}

Let $A$ be an algebra and $G$ be a group with identity $\sigma_1$.
We consider an \textbf{action} of $G$ on $A$, that is a
function
$\xymatrix@1@C=15pt{G\times A \ar[r] & A}$,
$\xymatrix@1@C=15pt{(\sigma,a) \ar@{|->}[r] & \sigma(a)}$,
such that:
\begin{enumerate}
\item[(a)] For each $\sigma$ in $G$, the map
$\sigma:\xymatrix@1@C=15pt{A \ar[r] & A}$
is an automorphism of algebra;
\item[(b)] $(\sigma \sigma')(a)=\sigma(\sigma'(a))$
for all $\sigma,\sigma'\in G$ and $a\in A$;
\item[(c)] $\sigma_1(a)=a$ for all $a\in A$.
\end{enumerate}

For any such action, the \textbf{skew group algebra} $A[G]$ is the
free left $A$-module with basis all the elements in $G$ endowed
with the multiplication given by
$(a\sigma)(b\sigma')=a\sigma(b)\sigma\sigma'$ for all $a,b\in A$
and $\sigma,\sigma' \in G$. Clearly, $A[G]$ admits a structure of
right $A$-module. Observe that $A[G]$ is
generally not connected and basic, but this will not play any
major role in the sequel.

In addition, any action of $G$ on $A$ induces a group action on
$\mod A$: for any $M\in \mod A$ and $\sigma\in G$, let
$^\sigma M$ be the $A$-module with the additive structure of $M$
and with the multiplication $a\cdot m=\sigma^{-1}(a)m$, for $a\in
A$ and $m\in M$. Given a morphism of $A$-modules
$\mor{f}{M}{N}$, define
$\mor{^{\sigma}f}{{}^{\sigma}M}{{}^{\sigma}N}$ by
$^{\sigma}f(m)=f(m)$ for each $m\in {}^{\sigma}M$. This defines an action of $G$ on $\mod A$; see \cite{ALR07}.

%In order to avoid confusion, and since the $A$-modules $M$ and
%$^{\sigma}M$ have the same elements but different external
%multiplication, we shall use the symbol $"\cdot"$ to express the
%multiplication $a\cdot m$ of an element $a$ in $A$ by an element
%$m$ in $^{\sigma}M$, while no symbol will be used to express the
%multiplication of $a$ by $m$ when $m$ is viewed as an element in
%$M$. Hence, the expressions $am$ and $a\cdot m$ are totally
%different.

When $G=\{\sigma_1, \sigma_2, \dots, \sigma_n\}$ is a finite group, the natural inclusion of $A$ in $A[G]$ induces the change of ring
functors
$F=\R\otimes_A - :\xymatrix@1@C=15pt{\mod A \ar[r] & \mod\R}$
and
$H=\Hom_{\R}(\R, -) :\xymatrix@1@C=15pt{\mod \R \ar[r] & \mod A}$.
These have been extensively studied in \cite{ALR07,RR85,Smith07},
for instance, assuming the order of $G$ is not a multiple of the characteristic of $k$.
We recall the following facts from \cite[(1.1)(1.8)]{RR85}.

%*********** Remark : rem RR
\fait{my_rem}{rem RR}{\emph{%
\begin{enumerate}
\item[(a)] $(F,H)$ and $(H,F)$ are adjoint pairs of
functors.
\item[(b)] Let $M\in\mod A$ and $\sigma\in G$. The subset
$\sigma\otimes_A M =\{\sigma\otimes_A m \ | \ m\in M\}$ of $FM$
has a structure of $A$-module given by
$a(\sigma\otimes_A m) =
\sigma\sigma^{-1}(a)\otimes_A m = \sigma \otimes_A
\sigma^{-1}(a)m= \sigma \otimes_A (a\cdot m),$
so that $\sigma\otimes_A M$ and $^{\sigma}M$ are isomorphic as
$A$-modules. 
%In addition, any element $x\in FM=\F{M}$ can be
%written in the form $x=\sum_{i=1}^{n}\sigma_i\otimes_A
%x_{\sigma_i}$, and in a unique way since the $\sigma_i$ form a
%basis for $\R$ as $A$-module. 
Therefore, as $A$-modules, we have
$$FM\cong \bigoplus_{i=1}^{n}(\sigma_i\otimes_A M) \cong
\bigoplus_{i=1}^{n} {}^{\sigma_i}M.$$
Then, 
$HFM\cong \bigoplus_{i=1}^{n}(\sigma_i\otimes_A M) \cong
\bigoplus_{i=1}^{n} {}^{\sigma_i}M.$ 
\item[(c)] Given a morphism $\mor{f}{M}{N}$ and $\sigma\in G$,
the map $\mor{^{\sigma}f}{{}^{\sigma}M}{^{\sigma}N}$ becomes
$\xymatrix@1@C=15pt{\sigma\otimes M \ar[r]& \sigma\otimes N}$ :
$\map{\sigma\otimes m}{\sigma\otimes f(m)}$ (that we also denote
$^{\sigma}f$) under the isomorphism $^{\sigma}M\cong \sigma\otimes
M$. In the sequel, we will freely use both notations. 
\end{enumerate}}}

We also need the following key observation.

%********** Remark : rem endo
\fait{my_rem}{rem endo}{\emph{%
In what follows, it will be convenient to consider the case where
a group $G$ acts on a Krull-Schmidt category $\A$ (with, say,
$\mor{^{\sigma}(-)}{\A}{\A}$ for $\sigma\in G$) and $T$ is a
(basic) $G$-stable object in $\A$, that
is, for each $\sigma\in G$, we have an isomorphism $\mor{\alpha_{\sigma}}{T}{{}^{\sigma}T}$. In this case, $\End_{\A}T$ is naturally endowed with an action of $G$, given by $\sigma(f)=\alpha_{\sigma}^{-1}\circ {}^{\sigma}f \circ \alpha_{\sigma}$, for $f\in \End_{\A}T$.
}}

%It is worthwhile to observe that if a group $G$ acts on the
%additive category $\add T$, for some $A$-module $T$, and $T$ is
%\textbf{multiplicity-free}, that is if any two distinct
%indecomposable direct summand of $T$ are not isomorphic, then $T$
%is $G$-stable under the given action.  However, in general, the
%converse fails to be true.

%-------------------------------------------------------------------------------
%
%  SECTION 2 : GROUP ACTIONS AND $G$-COMPATIBLE DERIVED EQUIVALENCES
%
%-------------------------------------------------------------------------------
\section{Group actions and $G$-compatible derived equivalences}
    \label{Group actions and $G$-compatible derived equivalences}

In this section, we recall how an action of $G$ on an additive
category $\A$ induces an action of $G$ on the homotopy and derived
categories of $\A$. Once this is done, we show that the
equivalences of derived categories induced by $G$-stable tilting
modules are $G$-compatible.

%-------------------------------------------------------------------------------
%  SUBSECTION : GROUP ACTIONS ON HOMOTOPY AND DERIVED CATEGORIES
%-------------------------------------------------------------------------------
\subsection{Group actions on homotopy and derived categories}
    \label{Group actions on homotopy and derived categories}

Let $G$ be a group and assume that $\A$ is an additive category on
which $G$ acts. For each $\sigma\in G$, let
$^\sigma (-):\xymatrix@1@C=15pt{\A \ar[r] & \A}$
be the automorphism of $\A$ induced by $\sigma$.
For any complex $\M=(M^i,
d^i_\M)_{i\in\mathbb{Z}}$ over $\A$ and $\sigma\in G$, let
$^{\sigma}\M$ be the complex $(^{\sigma}M^i,
{}^{\sigma}d^i_\M)_{i\in\mathbb{Z}}$. Moreover, given another
complex $\N=(N^i, d^i_\N)_{i\in\mathbb{Z}}$ and a morphism of
complexes $f=(\mor{f^i}{M^i}{N^i})_{i\in\mathbb{Z}}$, let
$^{\sigma}f=(\mor{^{\sigma}f^i}{{}^{\sigma}M^i}{{}^{\sigma}N^i})_{i\in\mathbb{Z}}$.
Clearly, $^{\sigma}f$ is a morphism of complexes.
Since
$^\sigma (-):\xymatrix@1@C=15pt{\A \ar[r] & \A}$
is an automorphism, this construction is compatible with the
homotopy relation. This allows to define, for each
$\sigma\in G$, an endomorphism
$\mor{^{\sigma}(-)}{K^b(\A)}{K^b(\A)}$. Moreover, since this action preserves the quasi-isomorphisms, it extends to an action on $D^b(\A)$.

%************* Prop : prop DA
\fait{my_prop}{prop DA}{Let $\sigma\in G$. The mapping
$\xymatrix@1@C=15pt{\M\ar@{|->}[r]& ^{\sigma}\M}$
(where $\M$ is a complex over $\A$) induces an action of $G$ on
$D^b(\A)$. In addition, the automorphisms
$\mor{^{\sigma}(-)}{D^b(\A)}{D^b(\A)}$ induced by the elements
$\sigma\in G$ are triangle-equivalences.}

At this point, recall that if $\A=\mod A$, for some finite
dimensional $k$-algebra $A$ of finite global dimension (for
instance if $A$ is piecewise hereditary \cite[(1.2)]{HRS96II}),
then $D^b(\A)$ has almost split triangles.  We have the following
result.

%********** Proposition : prop tps
\fait{my_prop}{prop tps}{Let $\A$ be as above and $\sigma\in G$.
Then, the automorphism $\mor{^{\sigma}(-)}{D^b(\A)}{D^b(\A)}$
preserves the almost split triangles.}

We get the following corollary, where the proof follows directly from
(\ref{prop tps}).

%********* Corollary : cor orbits
\fait{my_cor}{cor orbits} {Let $\A$ be as above and $\sigma\in G$.
Then the Auslander-Reiten translation $\tau$ and the functor
$^{\sigma}(-)$ commute on objects. In
particular, the functor $^{\sigma}(-)$ preserves the $\tau$-orbits in the Auslander-Reiten quiver of $D^b(\A)$.}

%-------------------------------------------------------------------------------
%
%  SUBSECTION 3 : G-COMPATIBLE DERIVED-EQUIVALENCES
%
%-------------------------------------------------------------------------------
\subsection{$G$-compatible derived equivalences}
    \label{$G$-compatible derived equivalences}

It is well-known from \cite{H88} that any tilting module induces
an equivalence of derived categories. Here, we show that the $G$-stable tilting modules
induce $G$-compatible equivalences. We recall the following facts
from \cite[(III.2)]{H88} : let $A$ be a finite dimensional $k$-algebra of finite global
dimension. Given a tilting $A$-module $T$ and $B=\End_AT$, the
functors
\begin{enumerate}
\item[(i)] $\mor{\Hom_A(T,-)}{K^b(\add T)}{{K^b(\proj B)}}$
\item[(ii)] $\xymatrix@1@C=15pt{\rho: K^b(\proj B) \
\ar@{^{(}->}[r] & K^b(\mod B) \ar[r] & \DBB}$
\item[(iii)] $\xymatrix@1@C=15pt{\phi: K^b(\add T) \
\ar@{^{(}->}[r] & K^b(\mod A) \ar[r] & \DBA}$
\end{enumerate}
are equivalences of triangulated categories, and the
composition
\begin{enumerate}
\item[(iv)] $\mor{\mbox{RHom}(T, -)=\rho \circ \Hom_A(T,-) \circ
\phi^{-1}}{\DBA}{\DBB}$
takes $T$ to $B$.
\end{enumerate}

Observe that above, and below, the functor $\Hom_A(T,-)$ is the
component-wise functor taking a complex $\T=(T^i,$ $f^i)$
in $K^b(\add T)$ to the complex
$\Hom_A(T, \T) = (\Hom_A(T, T^i),$ $\Hom_A(T, f^i))$
in $K^b(\proj B)$. 

%********* Proposition : prop Happel
\fait{my_prop}{prop Happel}{Let $A$ be an algebra and $G$ be a
group acting on $A$.  If $T$ is a tilting $A$-module which is
$G$-stable with respect to the induced action of $G$ on $\mod A$,
then the equivalences (i)-(iv) given above are $G$-compatible.}
\begin{proof}
First, let $\mor{\theta}{G}{\Aut A}$ be an action of $G$ on $A$
and, for each $\sigma\in G$, let $\mor{^{\sigma}(-)}{\mod A}{\mod
A}$ be the induced automorphism. Moreover, let $T$ be a tilting
$A$-module which is $G$-stable with respect to this action,
endowed with isomorphisms $\mor{\alpha_\sigma}{T}{\sigma \otimes
T}$ for $\sigma\in G$ as in (\ref{rem endo}).
The $G$-stability of $T$ gives rise to a natural action of $G$ on
$\add T$.  Then, following Section \ref{Group actions on homotopy
and derived categories}, the additive category $K^b(\add T)$
inherits a (component-wise) action of $G$.  We also denote the
induced automorphism on $K^b(\add T)$ by ${^{\sigma}(-)}$, for
$\sigma\in G$. Also, since $T$ is $G$-stable, it
follows from (\ref{rem endo}) that $B=\End_AT$ is endowed with
a natural action of $G$, which we extend to $K^b(\proj B)$ and
$D^b(B)$. Again, we denote the induced automorphisms by
$^{\sigma}(-)$, for $\sigma\in G$.

\noindent (i) $\mor{\Hom_A(T,-)}{K^b(\add T)}{{K^b(\proj B)}}$.
Let $\sigma\in G$ and $\T=(T^i,$ $f^i)$
be a complex in $K^b(\add T)$. We need to verify that $\sigma
\otimes_B \Hom_A(T, \T)\cong \Hom_A(T,\sigma \otimes \T)$
functorially. To do so, for each $i$ consider the map $\beta^i$ from $\sigma \otimes_B \Hom_A(T,T^i)$ to $\Hom_A(\sigma\otimes_A T, \sigma\otimes_A
T^i)$ taking $\sigma \otimes g$ onto $^{\sigma}g$
where $^{\sigma}g$ is as in (\ref{rem RR})(c).
Then $\beta^i$ is a morphism of $B$-modules. Since
$\beta^i$ is clearly bijective, it is an isomorphism of
$B$-modules. Since the $\beta^{i}$'s commute with any morphism in $\add T$, we have a
functorial isomorphism of complexes $\beta^\bullet$.
Now, since $T$ is $G$-stable, we have $T\cong \sigma \otimes T$
and a functorial isomorphism
$$\gamma^\bullet: \xymatrix@1@C=15pt{\sigma\otimes_B
\Hom_A(T,\T)\ar[r]^-{\beta^\bullet} & \Hom_A(\sigma\otimes
T,\sigma\otimes \T) \ar[r]^-{\cong} & \Hom_A(T,\sigma\otimes
\T)},$$
showing that $\Hom_A(T,-)$ is $G$-compatible.

\noindent (ii) and (iii).  Since the inclusions and localization
functors are clearly $G$-compatible, so are $\rho$ and
$\phi$. \cqfd
\end{proof}

%-------------------------------------------------------------------------------
%
%  SECTION : PIECEWISE HEREDITARY ALGEBRAS REVISITED
%
%-------------------------------------------------------------------------------
\section{Piecewise hereditary algebras revisited}
    \label{Piecewise hereditary algebras revisited}

The aim of this section is to prove Theorem \ref{thm DLS2}.  Before doing so, we
need to recall some facts concerning skew group algebras and prove
preliminary results.

%-------------------------------------------------------------------------------
%  SUBSECTION : PRELIMINARY RESULTS
%-------------------------------------------------------------------------------
\subsection{Preliminary results}
    \label{Preliminary results}

Let $A$ be an algebra and $G=\{\sigma_1, \sigma_2, \dots,
\sigma_n\}$ be a finite group acting on $A$ whose order is not a
multiple of the characteristic of $k$. Let $\sigma_1$ be the unit
of $G$, $T$ be an $A$-module and $\mor{f}{FT}{FT}$ be a $k$-linear
morphism.  Since $FT\cong \bigoplus_{i=1}^{n}(\sigma_i\otimes_A
T)$, $f$ is given by a matrix $f=(f_{\sigma_i, \sigma_j})_{1\leq i,j\leq
n}$
where each $f_{\sigma_i, \sigma_j}$ is a morphism from
$\sigma_j\otimes T$ to $\sigma_i\otimes T$.  Clearly, $f$ is
$A$-linear if and only if each $f_{\sigma_i, \sigma_j}$ is
$A$-linear.  Now, if $f$ is $A$-linear, then it is $A[G]$-linear
if and only if $f(\sigma_k\sigma_j\otimes
t)=\sigma_kf(\sigma_j\otimes t)$ for all $\sigma_j,\sigma_k\in G$,
and quick computations show that it is the case if and only if
$f_{\sigma_k\sigma_i, \sigma_k\sigma_j}(\sigma_k\sigma_j\otimes
t)=\sigma_k\cdot f_{\sigma_i, \sigma_j}(\sigma_j\otimes t)$
for all $\sigma_i, \sigma_j, \sigma_k\in G$ and $t\in T$. In particular, $f$ is determined by  $\{f_{\sigma_1, \sigma_1}, \dots, f_{\sigma_1, \sigma_n}\}$.

%*********** Proposition : lem B
\fait{my_prop}{lem B}{Let $A$ and $G$ be as above, and $T$ be a
$G$-stable (basic) tilting $A$-module with isomorphisms
$\mor{\alpha_\sigma}{T}{\sigma\otimes T}$ as in (\ref{rem endo}).
Then,
\begin{enumerate}
\item[\emph{(a)}] $FT$ is a $\R$-tilting module.
\item[\emph{(b)}] $\End_{\R}FT\cong (\B)[G]$.
\end{enumerate}}
\begin{proof}
\noindent (a). Since $F$ is exact and
preserves the projectives by (\ref{rem RR}), $FT$ has projective dimension at most one. In
addition, we have
$$\Ext_{\R}^1(FT,FT)\cong D\Hom_{\R}(FT, \tau(FT))\cong
D\Hom_{\R}(FT, F(\tau T)),$$
by the Auslander-Reiten formula and \cite[(4.2)]{RR85}. By
adjunction, this latter group is nonzero if and only if $\Hom_A(T,
H{F{(\tau T)}})\cong \Hom_A(T, \oplus_{\sigma\in
G}{}^{\sigma}(\tau T))$ is nonzero, where the isomorphism follows
from (\ref{rem RR}). However, if $\mor{f}{T}{{}^{\sigma}(\tau T)}$
is a nonzero morphism, with $\sigma\in G$, then
$\mor{{}^{\sigma^{-1}}f\circ \alpha_{\sigma^{-1}}}{T}{\tau T}$ is
nonzero, a contradiction to $\Ext^{1}_{A}(T,T)=0$. So, $\Ext_{\R}^1(FT,FT)=0$.
Finally, any short exact sequence 
$\xymatrix@1@C=15pt{0\ar[r]& A\ar[r] & T_0 \ar[r] & T_1\ar[r] &
0}$
in $\mod A$, with $T_0, T_1\in\add T$ induces a short exact sequence
$\xymatrix@1@C=15pt{0\ar[r]& \R\ar[r] & FT_0 \ar[r] & FT_1\ar[r] &
0}$
in $\mod\R$.  So $FT$ is a tilting $\R$-module.

\noindent (b). In order to show that $\End_{\R}FT$ and $(\B)[G]$
are isomorphic algebras, we construct explicit inverse
isomorphisms between them.
Let $f\in \End_{\R}FT$, and assume that $f$ is given by a matrix
$f=(f_{\sigma_i, \sigma_j})_{1\leq i,j\leq n}$, where each
$f_{\sigma_i, \sigma_j}$ is a morphism from $\sigma_j\otimes T$ to
$\sigma_i\otimes T$.  For each $i$, let $\mor{f_i=f_{\sigma_1,
\sigma_i}\circ \alpha_{\sigma_i}}{T}{\sigma_1\otimes T=T}$, and
define
$\nu: \End_{\R}FT \rightarrow (\B)[G]$ by $\nu(f)=
\sum_{i=1}^nf_i\sigma_i$
for each $f\in \End_{\R}FT$.

Conversely, let $\sum_{i=1}^nf_i\sigma_i \in (\B)[G]$, and
consider the family of morphisms $\{f_{\sigma_1, \sigma_1},\dots, f_{\sigma_1, \sigma_n}\}$, where
$f_{\sigma_1, \sigma_i}:= f_i\circ\alpha^{-1}_{\sigma_i}$ for each
$i$.  These morphisms determine an
$A[G]$-linear map $f$ from $FT$ to $FT$. Hence, define
$\mu: (\B)[G] \rightarrow \End_{\R}FT$ by $\mu(\sum_{i=1}^nf_i\sigma_i)=f$
for each $\sum_{i=1}^nf_i\sigma_i\in (\B)[G]$.

Clearly, $\nu$ and $\mu$ are inverse constructions, preserving
sums and units. It remains to show that $\nu$
preserves the product. For $f=(f_{\sigma_i, \sigma_j})_{1\leq i, j\leq n}$ and $g=(g_{\sigma_i, \sigma_j})_{1\leq i, j\leq n}$, let $\nu(f)=\sum_{i=1}^nf_i\sigma_i$ and
$\nu(g)=\sum_{j=1}^ng_j\sigma_j$. Then, 
$$\begin{array}{rl}
\nu(f)\cdot \nu(g)& =\sum_{i=1}^n \sum_{j=1}^n
f_i\sigma_i(g_j)\sigma_i\sigma_j\\
&=\sum_{i=1}^n \sum_{j=1}^n
f_i\circ(\alpha_{\sigma_i}^{-1}\circ{}^{\sigma_i}g_j\circ \alpha_{\sigma_i})\sigma_i\sigma_j\\
&=\sum_{i=1}^n \sum_{j=1}^n
(f_{\sigma_1, \sigma_i}\circ{}^{\sigma_i}g_{\sigma_1, \sigma_j}\circ ({}^{\sigma_i}\alpha_{\sigma_j}\circ\alpha_{\sigma_i}))\sigma_i\sigma_j\\
&=\sum_{k=1}^n (\sum_{i=1}^n
(f_{\sigma_1, \sigma_i}\circ g_{\sigma_i, \sigma_k}\circ\alpha_{\sigma_k}))\sigma_k\\
&=\nu(f\cdot g)
\end{array}
$$
 \cqfd
\end{proof}

%-------------------------------------------------------------------------------
%  SUBSECTION : PROOFS OF THEOREM 1 AND 2
%-------------------------------------------------------------------------------
\subsection{Proof of Theorem \ref{thm DLS2}}
    \label{Proofs of Theorem 1 and Theorem 2}

Let $\mc{A}$ be an
arbitrary abelian category, and $M\in\mc{A}$ satisfying
$\Ext^1(M,M)=0$ and $\Ext^2(M,N)=0$ for all $N\in \mc{A}$. Then,
the \textbf{right perpendicular category} $M^{\bot}$ is the full
subcategory of $\mc{A}$ containing the objects $N$ satisfying
$\Hom(M,N)=\Ext^1(M,N)=0$.  We define dually the left
perpendicular category $^{\bot}M$. It was shown in \cite{GL91}
that $M^{\bot}$ and $^{\bot}M$ are again abelian categories.

We can now proceed with the proof. The cases $\mathcal{H}=\mod H$ and $\mathcal{H}=\mbox{coh}\mathbb{X}$ are treated separately. Observe that our proof is an adaptation of the proofs of the main results in \cite{HRS88} and
\cite{HRS96II}, respectively, and to which we freely refer in the
course of the proof.

%\bigskip

%*********************************************************************
%    PROOF OF THEOREM 1
%*********************************************************************
\noindent {\bf Proof of Theorem \ref{thm DLS2} :} \\%
\noindent (a). Clearly, (iii) implies (ii), while (ii) implies (i)
by an easy induction using (\ref{prop Happel}) and its dual.

Now, let $\mathcal{H}=\mod H$, for some hereditary algebra $H$, and suppose that (i) holds. Assume, without loss of
generality, that $A$ and $H$ are basic and connected.  Let
$\mor{E}{\DBH}{\DBA}$ be a $G$-compatible equivalence; we shall identify the module categories
$\mod A$ and $\mod H$ with their images under the natural
embeddings into $\DBA$ and $\DBH$, respectively.

Let $\M$ be an object of $\DBH$ such that $E\M$ is isomorphic to
$A$ and assume, without loss of generality, that
$\M=M_0\oplus M_1[1]\oplus\cdots\oplus M_r[r]$
where the $M_i$ are $H$-modules and $M_0\neq 0\neq M_r$. Note that
$\Hom_H(M_i, M_j)=0$ if $i\neq j$, and $\Ext_H^1(M_i, M_j)=0$ if
$i+1\neq j$. Also, since $A$ is $G$-stable, then so is $\M$
and thus, since the $M_i[i]$ lie in different degrees, each of
them is also $G$-stable.

We prove our claim by induction on $r$.  If $r=0$, then one can
check that $\M=M_0$ is a $G$-stable tilting $H$-module, and so $A\cong \End_{H}M_0$.
In addition, since the isomorphism is given by $E$, the actions of
$G$ on $A$ and $\End_H M_0$ coincide. Also, since $H$ is
hereditary, then $\M$ is splitting and the result follows.

Now, assume inductively that the result holds true whenever $r$
takes a smaller value, or $r$ takes the same value and $M_0$ has
less indecomposable direct summands.
We shall find it convenient to construct a sequence of separating
tilting modules instead of splitting tilting modules. We recall
that $T_i$ is a separating tilting $A_i$-module with endomorphism
ring $A_{i+i}$ if and only if $T_i$ is a splitting tilting
$A_{i+i}^{op}$-module with endomorphism ring $A_{i}^{op}$. 
%
%Since $A_i$ and $A_{i+1}$ have equivalent derived categories if and only
%if so do $A_i^{op}$ and $A_{i+1}^{op}$, it will be clear that this
%will not lead into any fallacies.

Let $L=M_1\oplus M_2\oplus\cdots\oplus M_r$.  Then, by
\cite[(Proposition 3)]{HRS88}, the subcategory $L^{\bot}$ of $\mod
H$ is equivalent to $\mod \Lambda$, for some finite dimensional
hereditary algebra $\Lambda$.  Moreover, the inclusion functor
$\xymatrix@1@C=15pt{\mod \Lambda \ar@{^{(}->}[r] & \mod H}$
is full, faithful and exact.  In what follows, we identify
$L^{\bot}$ and $\mod\Lambda$. Observe that $M_0\in\mod \Lambda$,
and in fact is a tilting $\Lambda$-module. Also, $\mod \Lambda$ is
$G$-stable since so is $L$.

Let $Q=\Hom_k(\Lambda, k)$, a minimal injective cogenerator for
$\mod \Lambda$, and let
$$\T=Q\oplus  M_1[1]\oplus M_2[2]\oplus\cdots\oplus M_r[r].$$
We observe that $Q$ is a $G$-stable $\Lambda$-module.  Indeed, it
is easily verified that if $I$ is an injective $\Lambda$-module,
then so is $^{\sigma}I$ for each $\sigma\in G$. Because each $
M_i$ is also $G$-stable, then so is $\T $.
Following \cite{HRS88}, $E\T $ is (isomorphic to) a
separating tilting $A$-module.  Moreover, $E\T $ is $G$-stable
since $E$ is $G$-compatible and $\T $ is $G$-stable.  Now let
$B=\End_A E\T $. By (\ref{prop Happel}), we have an equivalence
%of triangulated categories
%
$$E'': \xymatrix@1@C=15pt@R=0pt{\DBH \ar[r]^E & \DBA \ar[r]^{E'} & D^b(B)%
\\
\qquad  \T  \qquad \ar@{|->}[r] & \qquad E\T  \qquad \ar@{|->}[r]
& \qquad B \qquad}$$
which is $G$-compatible since so are $E$ and $E'$.
Let $Q_0$ be an simple indecomposable direct summand of $Q$ : such
a module exists since $\Lambda$ is hereditary.  Then, for each
$\sigma\in G$, the $\Lambda$-module $^{\sigma}Q_0$ is also simple
injective and $E''({}^{\sigma}Q_0)={}^{\sigma}E''(Q_0)$ is a
simple projective $B$-module.  Since $B$ is connected,
$E''({}^{\sigma}Q_0)$ is not an injective $B$-module unless we are
in the trivial case of a simple algebra. Observe moreover that
since each $^{\sigma}Q_0$ is a simple injective $\Lambda$-module,
then the set $\{{}^{\sigma}Q_0 \ | \ \sigma\in G\}$ is finite, and
we denote its cardinality by $n$.

Now, imitating the arguments of \cite{HRS88}, we find, for each
$\sigma\in G$, an object $R_\sigma$ in $\DBH$ isomorphic to
$U_\sigma[1]$, for some $H$-module $U_\sigma$, such that
$$E''R_\sigma \cong \tau^{-1}E''({}^{\sigma}Q_0)\cong
\tau^{-1}({}^{\sigma}(E''Q_0))\cong
{}^{\sigma}(\tau^{-1}E''Q_0).$$
So $\oplus_{\sigma\in G}E''R_\sigma$ is a $G$-stable $B$-module.
We let
$S=(\oplus_{\sigma\in G}E''R_\sigma)\oplus E''N,$
where $(\oplus_{\sigma\in G}{}^{\sigma}Q_0) \oplus N=\T $.

At this point, it is worthwhile to observe that
$(\oplus_{\sigma\in G}R_\sigma)\oplus N = N_0\oplus N_1[1]\oplus
M_2[2]\oplus\cdots\oplus M_r[r]$
for some $N_0, N_1\in \mod H$, and where, by definition of the
$R_\sigma$, $N_0$ has $n$ less indecomposable direct summands than
$M_0$. In addition, by construction, $S$ is a "generalized" APR-tilting $B$-module; what is important for our purpose, and easy to verify, is that $S$ is a separating tilting $B$-module.

%We shall show that $S$ is a separating tilting $B$-module.
%
%If $I$ is an injective $B$-module, then we have
%
%$$\Hom_B(I,\tau S)=\Hom_B(I,\oplus_{\sigma\in
%G}{}^{\sigma}E''Q_0)=0$$
%%
%since ${}^{\sigma}E''Q_0$ is simple projective for each $\sigma\in
%G$. By the Auslander-Reiten formula, $S$ as projective dimension
%at most one. Moreover, $\Ext_B^1(S,S)\cong D\Hom_B(S, \tau
%S)=\Hom_B(S,\oplus_{\sigma\in G}{}^{\sigma}E''Q_0)=0$.
%%
%Finally, the number of indecomposable direct summands in an
%indecomposable decomposition of $S$ is equal to that of $\T $,
%which is equal to that of $B$.  So $S$ is a tilting $B$-module.
%
%Moreover, $S$ is a separating tilting module since the only
%indecomposable $B$-modules lying in $\mc{F}(S)$, the torsion-free
%class induced by $S$, are the $E''{}^{\sigma}Q_0$, whereas
%$\mc{T}(S)$, the corresponding torsion class, is the additive
%subcategory generated by all remaining indecomposables. Indeed, if
%$V$ is an indecomposable $B$-module, then $V\notin \mc{T}(S)$ if and
%only if $\Ext_B^1(S,V)\cong D\Hom_B(V,\oplus_{\sigma\in
%G}{}^{\sigma}E''Q_0)$ does not vanish, and this occurs if and only
%if $V\cong {}^{\sigma}E''Q_0$ for some $\sigma\in G$.

Let $C=\End_B(S)$. By (\ref{prop Happel}), we have an equivalence
of triangulated categories
$$E''': \xymatrix@1@C=15pt@R=0pt{\DBH \ar[r]^{E''} &
D^b(B) \ar[r]^{\cong} & D^b(C)%
\\
(\oplus_{\sigma\in G}R_\sigma)\oplus N \ \ar@{|->}[r] & \qquad S
\qquad \ar@{|->}[r] & \qquad C \qquad}$$
which is $G$-compatible.  Moreover, as observed earlier, we have
$(\oplus_{\sigma\in G}R_\sigma)\oplus N = N_0\oplus N_1[1]\oplus
M_2[2]\oplus\cdots\oplus M_r[r]$
where $N_0=0$ or contains less indecomposable direct summands than
$M_0$. By induction hypothesis, $C$ is piecewise hereditary of type
$\mod H$ and, using the separating tilting modules $E\T $ and $S$,
and keeping in mind our preceding discussion on separating tilting
modules, so is $A$. This shows the equivalence of conditions (i),
(ii) and (iii), when $\mathcal{H}=\mod H$.

We now assume that $\mathcal{H}=\mbox{coh}\mathbb{X}$, for some weighted projective line $\mathbb{X}$, and that condition (i) holds.  We will give all details until the case $\mod H$ carries over. 

Let
$\mor{E}{\DDBH}{\DBA}$ be a $G$-compatible equivalence of
triangulated categories.
Let $\M$ be an object of $\DDBH$ such that $E\M$ is isomorphic to
$A$ and assume, without loss of generality, that
$\M= M_0\oplus M_1[1]\oplus\cdots\oplus M_r[r]$
for some $ M_i\in \mc{H}$, with $ M_0\neq 0\neq  M_r$. Note that
$\Hom_{\mc{H}}( M_i,  M_j)=0$ if $i\neq j$, and $\Ext_{\mc{H}}^1(
M_i,  M_j)=0$ if $i+1\neq j$.
Also, since $A$ is $G$-stable, so is $\M$, and thus,
since the $ M_i[i]$ lie in different degrees, each of them is
also $G$-stable.

We prove our claim by induction on $r$.  If $r=0$, then one can
check that $\M= M_0$ is a $G$-stable tilting object in $\mc{H}$,
and
$A\cong 
\End_{\mc{H}}  M_0$.
So $A$ is quasitilted. In addition, since the isomorphism is given
by $E$, the actions of $G$ on $A$ and $\End_{\mc{H}}  M_0$
coincide.

Assume inductively that the result holds true in all cases
where either $r$ takes a smaller value, or $r$ takes the same
value and either $ M_0$ or $ M_r$ has less indecomposable direct
summands.
Then, $M_0$ is a $G$-stable tilting object in the abelian category
$(\oplus_{i=1}^r M_i[i])^{\bot}$ and $M_r$ is a $G$-stable cotilting
object in the abelian category ${}^{\bot}(\oplus_{i=0}^{r-1}
M_i[i])$.  By \cite{HRS96II}, one of these categories is a
module category over a hereditary artin algebra $H$.  The situation
is then reduced to the case $\mc{H}=\mod H$, and we are done. This proves (a).

(b). Now, assume that the equivalent conditions of (a)
are satisfied and that $G$ is a finite group whose order is not a
multiple of the characteristic of $k$.  Then there exists a
tilting object $T$ in $\mc{H}$ and a sequence of algebras
$\End_{\mc{H}}T=A_0, A_1, \dots, A_n=A$ on which $G$ acts and a
sequence $T_0, T_1, \dots, T_{n-1}$ where $T_i$ is a $G$-stable
tilting or cotilting $A_i$-module with endomorphism ring
isomorphic to $A_{i+1}$ for each $i$. Moreover, by (\ref{lem B})
and its dual, $A_i[G]\otimes_{A_i}T_i$ is a tilting or cotilting
$A_i[G]$-module for each $i$ and
$\End_{A_i[G]}(A_i[G]\otimes_{A_i}T_i)\cong A_{i+1}[G]$. Now,
since the order of $G$ is not a multiple of the characteristic of
$k$, $A_0[G]$ is hereditary if $\mc{H}=\mod H$ by \cite[(1.3)]{RR85}, and quasitilted if $\mc{H}=\mbox{coh}\mathbb{X}$ by \cite[(III.1.6)]{HRS96}.  The statement thus follows from \cite{HRS88} and
\cite{HRS96II}, respectively. \cqfd

%-------------------------------------------------------------------------------
%  SECTION : PROOF OF THEOREM 3
%-------------------------------------------------------------------------------
\section{Main result}
    \label{Proof of Theorem 3}

\subsection{Proof of Theorem \ref{thm DLS3}}

Assume that the hypotheses of Theorem \ref{thm DLS3} are satisfied. In view of Theorem \ref{thm DLS2}, it
would be sufficient to prove Theorem \ref{thm DLS3} to show that 
%if $A$ is piecewise hereditary of
%type $\mc{H}$, for some Ext-finite hereditary abelian category
%with tilting objects $\mc{H}$, and $G$ is a finite group acting on
%$A$ and whose order is not a multiple of the characteristic of
%$k$, then 
there is a $G$-compatible equivalence between
$D^b(A)$ and $D^b(\mc{H})$. 
%If this will be shown to be true in
%the case where $\mc{H}=\cohX$, for some weighted projective line
%$\mathbb{X}$, this unfortunately fails when $\mc{H}=\mod H$ for
%some hereditary algebra $H$.  
In this section, we show that this holds for $\mc{H}=\cohX$, and show that when
$\mc{H}=\mod H$, it is however possible to construct a derived equivalent hereditary algebra $H'$ on which $G$ acts and for which there is a $G$-compatible equivalence between $\DBA$ and
$D^b(H')$.

The first situation we consider is that of a piecewise hereditary
algebra $A$ of type $\mc{H}=\mod H$, for some hereditary algebra
$H$.  Here, it will be sufficient to assume that $G$ is a torsion
group acting on $A$. Let $Q=(Q_0, Q_1)$ be a finite and
acyclic quiver such that $H\cong kQ$, where $Q_0$ and $Q_1$
respectively denotes the set of vertices and arrows of $Q$.
All directed components are isomorphic to $\mathbb{Z}Q$ as translation quivers, and so have
only finitely many $\tau$-orbits.
Our first aim is to show that any such directed component $\Ga$
admits a section which is stable under the induced action of $G$
on $\DBA$. Recall that a full and connected subquiver $\Omega$ of
$\Ga$ is a \textbf{section} if it contains no oriented cycles,
it intersects each $\tau$-orbit of $\Ga$ exactly
once and it is convex.

%********* Lemma : lem directed
By the above, we remark that, if we set $^{\sigma}\Ga :=
\{{}^{\sigma}\M \ | \ \M\in \Ga\}$,  then $^{\sigma}\Ga=\Ga$ for every $\sigma\in G$.
%We verify that by applying repeatedly the automorphism
%$^{\sigma}(-)$ yields to $\Ga_i=\Ga_{i+km_\sigma}$,
%where $\{\Ga_i\}_{i\in\mathbb{Z}}$ is the family of directed
%components of $\Gamma(D^b(A))$ for each
%$i\in\mathbb{Z}$. So, $km_\sigma=0$, and thus $m_\sigma=0$.
%\medskip
%
We now construct a $G$-stable section in $\Ga$ as follows.

%************** Definition : defn section
\fait{my_def}{defn section}{\emph{%
Let $G$, $A$ and $\Ga$ be as above, and let $\X$ be a fixed object
in $\Ga$.  We define $\Sigma$ $(= \Sigma_{\X})$ to be the full
subquiver of $\Ga$ formed by the objects $\M$ in $\Ga$ such that
there exists a path
$\xymatrix@1@C=15pt{^{\sigma}\X \ar@{~>}[r] & \M}$
for some $\sigma\in G$ and any such path is sectional.}}

%Observe that $\Sigma\subseteq\Ga$ by (\ref{lem directed})(a). We
%prove that $\Sigma$ is a section in $\Ga$.

%********* Lemma : lem orbits
\fait{my_lem}{lem orbits}{Let $\Sigma$ and $\Ga$ be as above. Then
$\Sigma$ intersects each $\tau$-orbit of $\Ga$ exactly once.}
\begin{proof}
Let $\M\in\Ga$.  Since $\Gamma$ is directed, for each $\sigma\in G$, there exists an integer $r_\sigma$ such that
there exists a path from $^{\sigma}\X$ to $\tau^{r}\M$ in $\Ga$ if
and only if $r\leq r_\sigma$.  Clearly, any path from
$^{\sigma}\X$ to $\tau^{r_\sigma}\M$ is sectional.  There exists an integer $s$ which is maximal for
the property that there exists a path
$\xymatrix@C=15pt{^{\sigma}\X \ar@{~>}[r] & \tau^s\M}$
in $\Ga$, for some $\sigma\in G$. The maximality of $s$ gives
$\tau^s\M\in\Sigma$.  The uniqueness of $\tau^s\M$ follows from
the definition of $\Sigma$. \cqfd
\end{proof}

%********* Lemma : lem convex
\fait{my_lem}{lem convex}{Let $\Sigma$ and $\Ga$ be as above.
Let
$\xymatrix@C=15pt{\omega : \M \ar@{-}[r] & M_1^{\bullet}
\ar@{-}[r] & \cdots \ar@{-}[r] & M_n^{\bullet}}$
be a walk in $\Ga$, with $n\geq 1$ and $\M\in \Sigma$. Then $\tau^k
M_n^{\bullet} \in \Sigma$ for some integer $k$. Moreover, $\M$ and
$\tau^k M_n^{\bullet}$ belong to the same connected component of
$\Sigma$.}
\begin{proof}
Let $\omega$ be as in the statement. We prove our claim by
induction.  First observe that by Lemma \ref{lem orbits}, it follows that if $\mor{f}{\M}{\N}$ is an irreducible morphism in $\Ga$, with $\M\in \Sigma$, then $\N\in\Sigma$ or $\tau\N\in\Sigma$, and dually. Hence, if $n=1$, then the claim follows from fullness of
$\Sigma$. Now, assume that the statement holds for
$n-1$. There exists $k\in\mathbb{Z}$ such that $\tau^k
M_{n-1}^{\bullet}$ belongs to the same connected component of
$\Sigma$ as $\M$. By translation, there exists an irreducible
morphism between $\tau^kM_{n-1}^{\bullet}$ and $\tau^k
M_n^{\bullet}$. Another application of the case $n=1$ gives the result.
\cqfd
\end{proof}

%********** Proposition : prop cs
\fait{my_prop}{prop cs}{The subquiver
$\Sigma$ is a $G$-stable section in $\Ga$.}
\begin{proof}
%Let $\Ga$ be a directed component of $\Gamma(\DBA)$, and
%$\Sigma\subseteq \Ga$ be defined as above.  We show that $\Sigma$
%is a section in $\Ga$.
First, $\Sigma$ is a full subquiver of $\Ga$ by definition.
Moreover, $\Sigma$ contains no oriented cycles (since $\Ga$ is
directed) and intersects each $\tau$-orbit of $\Ga$ exactly once by
(\ref{lem orbits}).  Since $\Sigma$ is clearly convex, it remains to show that $\Sigma$ is
connected and $G$-stable.
For the connectedness, assume that $\M,\N$ are two objects in
$\Sigma$. Since $\Ga$ is connected, there exists a walk from $\M$ to
$\N$ in $\Ga$. By (\ref{lem convex}), there exists $r\in\mathbb{Z}$
such that $\tau^r\N$ belongs to the same connected component of
$\Sigma$ as $\M$. Since $\Sigma$ intersects each $\tau$-orbit
exactly once, we get $r=0$, and so $\Sigma$ is connected.
Finally, $\Sigma$ is $G$-stable since, for each $\sigma\in G$, the
functor $\mor{^{\sigma}(-)}{\DBA}{\DBA}$ commutes with the
Auslander-Reiten translation $\tau$ by (\ref{cor orbits}), and
thus preserves the sectional paths. \cqfd
\end{proof}

%********* Remark : rem H-modules
%\fait{my_rem}{rem modules}{\emph{%
%Applying the suspension functor and the Auslander-Reiten
%translation if necessary, it is easily seen that, up to
%equivalence, the section $\Sigma$ can be taken in $(\mod H)[0]$,
%that is the set $\{M[0] \ | \ M \in \mod H\}$, which we identify
%with $\mod H$.}}

%We now mention that if a group $G$ acts on $\DBA$, for some
%algebra $A$, and there exists an equivalence of triangulated
%categories $\mor{E}{\DBA}{\DBB}$, for some algebra $B$, then one
%can carry the action of $G$ on $\DBB$ as follows : given an object
%$\N$ in $D^b(B)$ and $\sigma\in G$, set
%$^{\sigma}\N=E({}^{\sigma}(E^{-1}\N))$, and define the action on
%the morphisms in the same way.
%%
%The functor $E$ is then $G$-compatible since we have
%%
%$^{\sigma}(E(-))=E({}^{\sigma}(E^{-1}E(-)))=E({}^{\sigma}(-))$
%%
%for each $\sigma\in G$. With this construction in mind, we get the
%following key result.

%The following proposition shows how one can convert any
%triangle-equivalence between $\DBA$ and $D^b(H)$ into a
%$G$-compatible triangle-equivalence between $\DBA$ and $D^b(H')$,
%for a suitably chosen hereditary algebra $H'$ arising from a
%$G$-stable section.

%********* Proposition : prop hered
\fait{my_prop}{prop hered}{Let $A$ be a piecewise hereditary
algebra of type $\mod H$, for some hereditary algebra $H$, and $G$
be a torsion group.  Then, for any action of $G$ on $A$, there
exists a hereditary algebra $H'$ and an action of $G$ on $H'$
inducing a $G$-compatible equivalence of triangulated categories
between $\DBA$ and $D^b(H')$.}
\begin{proof}
Since $G$ is a torsion group, it
follows from (\ref{prop cs}) that $\DBA$ admits a $G$-stable section $\Sigma$.  Let $H'=\End_{D^b(A)}\Sigma$.  By Rickard's
Theorem \cite{R89}, there exists an equivalence of triangulated
categories
$E: \xymatrix@1@C=15pt@R=0pt{\DBA \ar[r] & D^b(H')}$
which takes $\Sigma$ to the full subquiver $\Omega$ of projective $H'$-modules in $D^b(H')$.
Under the identification $\DBA\cong D^b(H')$, $D^b(H')$ is
endowed with an action of $G$ and, for any $\sigma\in G$, we let
$\mor{_{\sigma}(-)}{D^b(H')}{D^b(H')}$ denote the induced
automorphism.  These automorphisms restrict to automorphisms
of $(\mod H')[i]$, for any $i\in\mathbb{Z}$, by
\cite[(IV.5.1)]{H88}. To prove our claim, it then remains
to show that there exists an action of $G$ on $H'$ such that the
induced action on $D^b(H')$ (see Section \ref{Group actions on homotopy and derived categories}) coincides with the
action carried from $D^b(A)$.

For this sake, observe that since $\Sigma$ is $G$-stable, so
is $\Omega$. Moreover, $\Omega$ is the ordinary quiver associated
to $H'$, and so $H'\cong k\Omega$. We define an action of $G$
on $H'$ as follows: let $\{e_1, e_2, \dots, e_n\}$ be a complete
set of primitive orthogonal idempotents of $H'$ and let $\{P_1,
P_2, \dots, P_n\}$ be the associated indecomposable projective
$H'$-modules, each of them being a vertex of $\Omega$. Then, for
$\sigma\in G$, set $\sigma(e_i)=e_j$ if $_\sigma P_i =P_j$.
Moreover, if $\alpha$ is an arrow of $\Omega$, then set
$\sigma(\alpha)={}_{\sigma}\alpha$. This defines an action of $G$
on $H'$, and further on $D^b(H')$. For each $\sigma\in G$, we let
$\mor{^{\sigma}(-)}{D^b(H')}{D^b(H')}$ denote the induced
automorphism. The equivalences
$_\sigma(-)$ and $^\sigma(-)$ coincide, up to a functorial
isomorphism, because they clearly coincide on projectives. \cqfd
\end{proof}

As we will see, the above proposition will play a major role in
the proof of Theorem \ref{thm DLS3}.
We now consider the case where  $A$ is a piecewise hereditary
algebra of type $\mc{H}=\cohX$, for some weighted projective line
$\mathbb{X}$, and $G$ is a group acting on $A$.
For more details concerning the categories of coherent sheaves on
a weighted projective lines, we refer to
\cite{GL87,LM00}.

Let $p_1, p_2, \dots, p_r$ be a set of natural numbers and
$\mathbb{X}=\mathbb{X}(p_1, p_2, \dots, p_r)$ be a weighted
projective line over $k$ of type $p_1, p_2, \dots, p_r$ (in the
sense of \cite{GL87}). Let $\mc{H}=\cohX$ be the category of
coherent sheaves on $\mathbb{X}$.  Then $\mc{H}$ is a hereditary
abelian category with tilting objects. It is known that
there exists a tilting object $T\in\mc{H}$ such that
$\End_{\mc{H}}T=C(p_1, p_2, \dots, p_r)$, where $C(p_1, p_2, \dots,
p_r)$ is a canonical algebra of type $p_1, p_2, \dots, p_r$ (in
the sense of \cite{R84}).
An important classification tool is the slope
function $\mor{\mu}{\mc{H}}{\mathbb{Q}\cup\infty}$; see \cite{LM00}.

Then, we get the following proposition, whose proof easily follows
from \cite[(4.4)]{LM00}.  We include a sketch of proof for the convenience
of the reader.

%*********** Prop : prop cohX
\fait{my_prop}{prop cohX}{Let $A$ be a piecewise hereditary
algebra of type $\emph{coh}\mathbb{X}$, for some weighted
projective line $\mathbb{X}$, and $G$ be a group.  For any action
of $G$ on $A$, there exists an action of $G$ on
$\emph{coh}\mathbb{X}$ and a $G$-compatible equivalence of
triangulated categories between $\DBA$ and
$D^b(\emph{coh}\mathbb{X})$.}
\begin{proof}
Assume that $G$ acts on $A$, and let
$\mor{^{\sigma}(-)}{\DBA}{\DBA}$ be the induced isomorphism for
each $\sigma\in G$. Also, let $P_1, P_2, \dots, P_n$ be a complete
set of indecomposable projective $A$-modules (up to isomorphism).
Since $A$ and $\cohX$ are derived equivalent, it follows from
Rickard's Theorem \cite{R89} that there exists a tilting complex
$T=T_1\oplus T_2 \oplus \cdots \oplus T_n$ such that
$A\cong\End_{D^b(\cohX)}T$. Moreover, we may assume that the
equivalence sends each indecomposable direct summand $T_i$ of $T$
to $P_i$, for $i=1,2,\dots, n$.  With this equivalence, $G$ acts
on $D^b(\cohX)$ and, for each $\sigma\in G$, the induced
automorphism $\mor{^{\sigma}(-)}{D^b(\cohX)}{D^b(\cohX)}$ yields a
permutation of $T_1, T_2, \dots, T_n$ hence of their slopes.  Then, by
\cite[(4.1)]{LM00}, $^{\sigma}(-)=T^m\circ f_\sigma$, where $T$ is
the translation functor of $D^b(\cohX)$ and $f_\sigma$ is an
automorphism of $\cohX$. Now, since $^{\sigma}(-)$ permutes $T_1,
T_2, \dots, T_n$, we further deduce that $m=0$, and thus
$^{\sigma}(-)$ restricts to $\cohX$.  This shows that $G$ acts on
$\cohX$, hence the above equivalence between $\DBA$ and
$D^b(\cohX)$ is $G$-compatible. \cqfd
\end{proof}

%*********************************************************************
%    PROOF OF THEOREM 3
%*********************************************************************
\noindent {\bf Proof of Theorem \ref{thm DLS3}} :
This follows from (\ref{prop hered}), (\ref{prop cohX}) and Theorem \ref{thm DLS2}.
\cqfd

%-------------------------------------------------------------------------------
%
%  SECTION 4.2 : EXAMPLES
%
%-------------------------------------------------------------------------------
\subsection{An example}
    \label{An example}

In this section, we illustrate Theorem \ref{thm DLS3} and the
mechanics of (\ref{lem B}) on a small example.
Let $A$ be the path algebra of the quiver $(1)$ below with relations
$\alpha\beta=0$ and $\alpha'\beta'=0$. The cyclic
group $G=\mathbb{Z}/2\mathbb{Z}$ acts on $A$ by switching $1$ and
$1'$, $3$ and $3'$, $\alpha$ and $\alpha'$, $\beta$ and $\beta'$,
and fixing the vertex $2$. By applying the method explained in
\cite[(Section 2.3)]{RR85}, we get that the skew group algebra
$A[G]$ is (Morita equivalent to) the path algebra of the quiver
$(2)$ below with relations $\gamma\delta=\gamma'\delta'$.
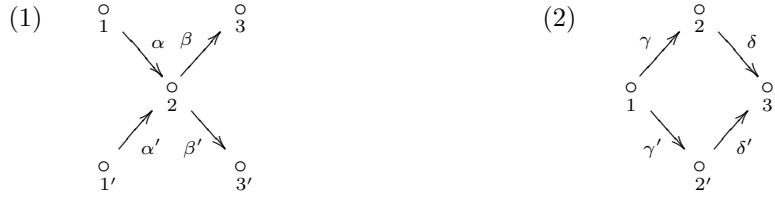
$$\xymatrix@R=12pt@C=12pt{(1)&\underset{1}{
\circ}\ar@<0.6ex>[dr]^{\alpha} && \underset{3}{
\circ} &&&&&&(2)& &\underset{2}{\circ} \ar@<0.6ex>[dr]^{\delta}&\\ %
&& \underset{2}{\circ} \ar@<0.6ex>[ur]^{\beta} \ar@<0.6ex>[dr]_{\beta'}&&&&&&&&\underset{1}{\circ} \ar@<0.6ex>[ur]^{\gamma} \ar@<0.6ex>[dr]_{\gamma'}  && \underset{3}{\circ}&\\ %
&\underset{\ 1'}{\circ} \ar@<0.6ex>[ur]_{\alpha'} && \underset{\
3'}{\circ}&&&&&&&& \underset{\ 2'}{\circ} \ar@<0.6ex>[ur]_{\delta'}&} $$
On the other hand, the Auslander-Reiten quiver of $\DBA$ consists
of a unique directed component $\Ga$ given as follows, where the
pair $(M,N)_n$ indicates that the homology in degree $n$ is $M$,
and the homology in degree $n+1$ is $N$, for some
$n\in\mathbb{Z}$. The $A$-modules $M$ and $N$ are represented by
their Loewy series.
$$ 
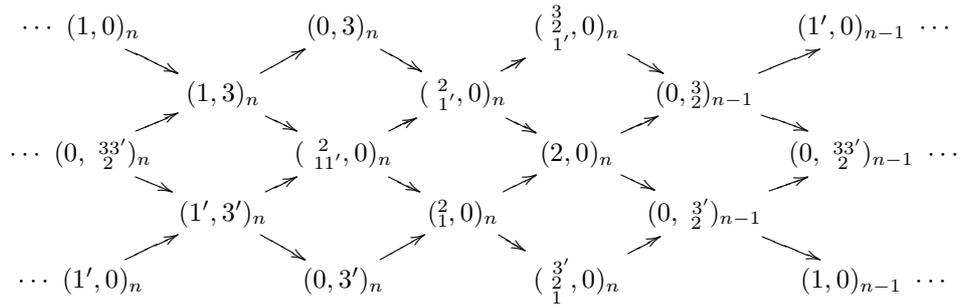
\xymatrix@R=5pt@C=4pt{%
\cdots \ (1,0)_n \ar[dr] && (0,3)_n\ar[dr] && (\substack{3
\\ 2\\ \ 1'}, 0)_n \ar[dr] &&
(1',0)_{n-1} \ \cdots \\
& (1,3)_n \ar[ur] \ar[dr] && (\substack{2 \\ \ 1'}, 0)_n \ar[ur]
\ar[dr] && (0,\substack{3 \\ 2})_{n-1} \ar[ur] \ar[dr] \\
\cdots \ (0, \substack{\ 33' \\ 2})_n \ar[ur] \ar[dr] &&
(\substack{2
\\ \ 11'}, 0)_n \ar[ur] \ar[dr] && (2,0)_n \ar[ur] \ar[dr] && (0,\substack{ \ 33' \\ 2})_{n-1} \ \cdots \\
& (1',3')_n \ar[ur] \ar[dr] && (\substack{2 \\ 1}, 0)_n \ar[ur]
\ar[dr] && (0,\substack{\ 3' \\ 2})_{n-1} \ar[ur] \ar[dr] \\
\cdots \ (1',0)_n \ar[ur] && (0,3')_n\ar[ur] && (\substack{\ 3'
\\ 2\\ 1}, 0)_n \ar[ur] &&
(1,0)_{n-1} \ \cdots}$$
Now, let $\X$ be a fixed object in the above directed component
$\Ga$, say $\X=(1, 3)_n$, and construct, as in (\ref{defn
section}), the unique section $\Sigma=\Sigma_{\X}$ of $\Ga$ having
the objects of the form $^{\sigma}\X$ as sources, for $\sigma\in
G$. Clearly, the induced action of $G=\mathbb{Z}/2\mathbb{Z}$ of
$\DBA$ switches the objects $(1,3)_n$ and $(1', 3')_n$, and so
$\Sigma$ is the full subquiver of $\Ga$ generated by the objects
$(1,3)_n, (1',3')_n, (0,3)_n, (0,3')_n$ and $(\substack{2 \\ \
11'}, 0)_n$. Now, let $H'=\End \Sigma$.

Then, by Theorem
\ref{thm DLS2}, there exist a sequence $H'=A_0, A_1, \dots, A_n=A$
of algebras and a sequence $T_0, T_1, \dots, T_{n-1}$ of modules
such that, for each $i$, $A_{i+1}=\End_{A_i}T_i$ and $T_i$ is a
$G$-stable tilting $A_i$-module.  Here, $n=1$, and so $A$ is a
tilted algebra of type $H'$. Indeed, the Auslander-Reiten quiver
of $H'$ is given by
$$ 
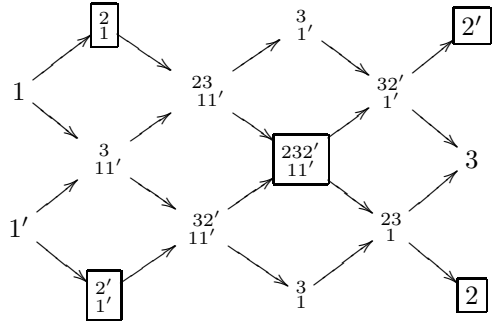
\xymatrix@R=7pt@C=15pt{%
& *+[F]{\substack{2\\ 1}} \ar[dr] && {\substack{3 \\ \ 1'}}
\ar[dr] && *+[F]{2'} \\
1 \ar[ur] \ar[dr] && {\substack{23 \\ \ \ 11'}} \ar[ur] \ar[dr] &&
{\substack{ 32' \\ 1'}} \ar[ur] \ar[dr]\\
& {\substack{3 \\ \ 11'}} \ar[ur] \ar[dr] &&
*+[F]{{\substack{232'\\11'}}} \ar[ur] \ar[dr] && 3 \\
1' \ar[ur] \ar[dr] && {\substack{ \ 32' \\ 11'}} \ar[ur] \ar[dr]
&& {\substack{23\\ 1}} \ar[ur] \ar[dr] \\
& *+[F]{\substack{2'\\ 1'}} \ar[ur] && {\substack{ 3 \\ 1}}
\ar[ur] && *+[F]{2}}$$
and it is easily seen that if $T$ is the direct sum of the identified
indecomposable modules in the above diagram, then $T$ is a
$G$-invariant tilting $H'$-module such that $\End_{H'}T\cong A$.
Now, since $G$ is cyclic, the method explained in \cite[(Section
2.3)]{RR85} gives that $H'[G]$ is (Morita equivalent to) the path
algebra of the quiver $(3)$ below, and the Auslander-Reiten quiver of $H'[G]$ is given by the quiver $(4)$ below.
$$
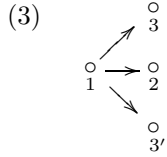
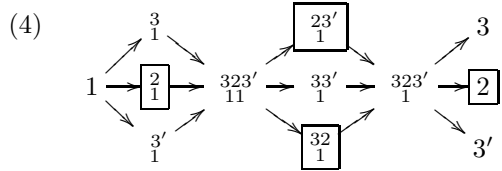
\xymatrix@R=5pt@C=10pt{%
(3)&& \underset{3}{\circ}&&&&&&&(4)&& {\substack{3\\ 1}} \ar[dr] && *+[F]{{\substack{\ 23' \\ 1}}}
\ar[dr] && 3\\%
& \underset{1}{\circ} \ar@<0.6ex>[ur] \ar@<0.6ex>[dr] \ar@<0.6ex>[r] & \underset{2}{\circ}&&&&&&&&1 \ar[r] \ar[ur] \ar[dr] & *+[F]{{\substack{2\\ 1}}} \ar[r] &
{\substack{\ 323' \\ 11}} \ar[r] \ar[ur] \ar[dr] & {\substack{\
33'\\ 1}} \ar[r] & {\substack{\ 323' \\ 1}} \ar[r] \ar[ur] \ar[dr]
& *+[F]{2}\\%
 && \underset{\ 3'}{\circ}&&&&&&&&& {\substack{\ 3'\\ 1}} \ar[ur] && *+[F]{{\substack{32
\\ 1}}} \ar[ur] && 3'}$$
where the identified indecomposable modules in the quiver $(4)$ correspond with the
indecomposable direct summand of the tilting $H'[G]$-module $FT$
of (\ref{lem B}). It is then easily verified, as predicted by
(\ref{lem B}), that $\End_{H'[G]}FT\cong (\End_{H'}T)[G]\cong
A[G]$.

%-------------------------------------------------------------------------------
%
%   ACKNOWLEDGEMENTS
%
%-------------------------------------------------------------------------------
%\bigskip
\noindent {\bf ACKNOWLEDGEMENTS.} The authors thank Ibrahim Assem
who suggested the problem considered in this paper. This work was
started while the second author was visiting the Université de
Sherbrooke and Bishop's University under a postdoctoral fellowship, and
pursued while the first and third authors were consecutively visiting the Universidad de la
Rep\'ublica del Uruguay. They would like to express
their gratitude to Marcelo for his warm hospitality. Finally, the
third author thanks Helmut Lenzing for a brief but fruitful
discussion while he was visiting the NTNU in Norway in January 2007.

%-------------------------------------------------------------------------------
%
%  BIBLIOGRAPHIE
%
%-------------------------------------------------------------------------------

\bibliography{biblioDLS}

\end{document}